%% file: main.tex
\documentclass{amsart}

\usepackage[T1]{fontenc}
\usepackage[utf8]{inputenc}
\usepackage[USenglish]{babel}
\usepackage{textcase}

\usepackage[
	bookmarks=true,
	plainpages=false,
	linktocpage,
	colorlinks=true,
	citecolor=green!80!black,
	linkcolor=red!70!black,
	filecolor=magenta,
	urlcolor=magenta,
	breaklinks,
	pdfauthor={Martin Winter},
]{hyperref}

\usepackage{amsmath,amsthm,amssymb}
\usepackage{calc, mathtools} %calc for \widthof, mathtools for \mathclap

\usepackage{colortbl,color,xcolor} % colors
\usepackage{centernot} % crossing out symbols correctly
\usepackage{array} % arrays
\usepackage{enumitem,moreenum} % numbered lists where each item can be labeled
\usepackage{cite} % for grouping references in the text which are adjacent in the bibliography
\usepackage[nameinlink,capitalize,noabbrev]{cleveref}
\usepackage{nicefrac}

\usepackage{tikz-cd} 

\usepackage[font=small,labelfont=bf]{caption}

\input{header}

\begin{document}

%%%%%%%%%%%%%%%%%%%%%%%%%%%%%%%%%%%%%%%%%%%%%%%%%%%%%%%%%%%%%%%%%%%%%%%%%%%%%%%%%%%%%%

\expandafter\title
{Classification of vertex-transitive zonotopes}
		
\author[M. Winter]{Martin Winter}
\address{Faculty of Mathematics, University of Technology, 09107 Chemnitz, Germany}
\email{martin.winter@mathematik.tu-chemnitz.de\newline\rule{0pt}{1.5cm}\includegraphics[scale=0.7]{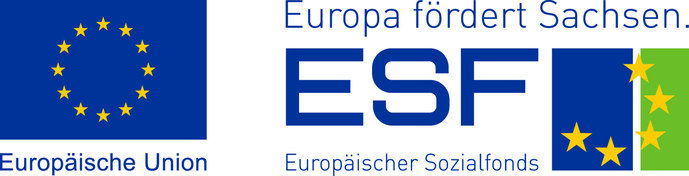}}
	
\subjclass[2010]{51F15, 52B11, 52B12, 52B15}
\keywords{convex polytopes, zonotopes, permutahedra, symmetry, vertex-transitivity, reflection groups, root systems}
		
\date{\today}
\begin{abstract}
%We investigate zonotopes with a high degree of regularity or symmetry in the structure of their $\delta$-faces for $\delta\in\{0,1,2\}$. 
%Here, regularity means that all the $\delta$-faces of the zonotope share common properties, while symmetry means that the symmetry group of the zonotope acts transitively on the $\delta$-faces.

We give a full classification of vertex-transitive zonotopes.

We prove that a vertex-transitive zonotope is a $\Gamma$-permutahedron for some finite reflection group $\Gamma\subset\Ortho(\RR^d)$.
The same holds true for zonotopes in which all vertices are on a common sphere, and all edges are of the same length (which we call \emph{homogeneous} zonotopes).
The classification of these then follows from the classification of finite reflection groups.

We proof that root systems can be characterized as those centrally symmetric sets of vectors, for which all intersections with half-spaces, that contain exactly half the vectors, are congruent.
We provide a further sufficient condition for a centrally symmetric set being a root system.

%This enables us to proof an alternative characterization of root systems.
%We further give a full classification of edge-transitive zonotopes, in particular, we show that in $\ge 4$ dimensions there are only $d$-cubes and $(2k,\...,2k)$-hyperprisms.
%Similar results are obtained if vertex-transitivity/edge-transitivity is replaced by a weaker form of regularity.
%
%We also show that zonotopes with congruent faces are linked to arrangements of equiangular lines. 
%(which seems to be known).

%We prove that all vertex-transitive zonotopes are $\Gamma$-permutahedra for some reflection group $\Gamma$, and that the same already holds if we only require to have vertices on a common sphere, and simultaneously, all edges of the same length. 
%We prove that besides $d$-cubes and some polygons, there are only two edge-transitive (actually edge-regular) zonotopes, and both are 3-dimensional. We show that the 2-face-transitive zonotopes in $\ge 3$ dimensions are in a one-to-one correspondence with arrangements of equiangular lines.
\end{abstract}

\maketitle

%\setcounter{tocdepth}{1}
%\tableofcontents

\input{\secdir introduction}

\input{\secdir basics}

\input{\secdir vertex-transitive}

\input{\secdir conclusion}

\par\medskip
\noindent
\textbf{Acknowledgements.} The author gratefully acknowledges the support by the funding of the European Union and the Free State of Saxony (ESF).

\bibliographystyle{unsrt}
\bibliography{literature}

\newpage
\input{\secdir appendix}

\end{document}

%% file: header.tex
\newcommand{\NN}{\mathbb{N}}    % natural numbers            
\newcommand{\RR}{\mathbb{R}}    % real numbers                      
\newcommand{\Sph}{\mathbb{S}}    % sphere 
                   
\newcommand{\F}{\mathcal{F}}    % faces
\newcommand{\Tsymb}{\top}
\newcommand{\T}{^{\Tsymb}}

\newcommand{\nozero}{\setminus\{0\}}

\def\^#1{^{(#1)}}
\def\s^#1{^{\smash{(#1)}}}

\newcommand{\wideword}[1]{\quad\text{#1}\quad}
\newcommand{\wideand}{\wideword{and}}

\def\:{\colon}

\newcommand{\cupdot}{\mathbin{\mathaccent\cdot\cup}}

\newcommand{\mylabel}{$(\roman*)$}
\newenvironment{myenumerate}{\begin{enumerate}[label=\mylabel]}{\end{enumerate}}

\newcommand{\freespace}{\kern.07em} 
\newcommand{\free}{\freespace\cdot\freespace} 

\newcommand{\enquote}[1]{``#1''}                                 % quotation marks

\newcommand{\msays}[1]{{\footnotesize\textcolor{red}{#1}}}
\newcommand{\TODO}{\msays{TODO}}

\theoremstyle{plain}  % theorem, lemma, corollary, proposition, conjecture, criterion, algorithm
\newtheorem{theorem}{Theorem}[section]
\newtheorem{corollary}[theorem]{Corollary}
\newtheorem{lemma}[theorem]{Lemma}
\newtheorem{proposition}[theorem]{Proposition}

\theoremstyle{definition} % definition, condition, problem, example
\newtheorem{definition}[theorem]{Definition}

\newtheorem{remark}[theorem]{Remark}
\newtheorem{question}[theorem]{Question}
\newtheorem{observation}[theorem]{Observation}

\crefname{theorem}{Theorem}{Theorems}
\crefname{proposition}{Proposition}{Propositions}
\crefname{lemma}{Lemma}{Lemmas}
\crefname{corollary}{Corollary}{Corollaries}
\crefname{remark}{Remark}{Remarks}
\crefname{example}{Example}{Examples}
\crefname{definition}{Definition}{Definitions}
\crefname{problem}{Problem}{Problems}
\crefname{observation}{Observation}{Observation}
\crefname{construction}{Construction}{Construction}

\DeclareMathOperator{\conv}{conv}

\DeclareMathOperator{\Aut}{Aut}

\DeclareMathOperator{\Ortho}{O}

\DeclareMathOperator{\Span}{span}

\DeclareMathOperator{\Id}{Id}

\DeclareMathOperator{\Star}{Star}
\DeclareMathOperator{\Gen}{Gen}
\DeclareMathOperator{\Zon}{Zon}
\DeclareMathOperator{\Orb}{Orb}
\DeclareMathOperator*{\argmax}{arg\,max}

\let\eps=\epsilon
\let\eset=\varnothing

\let\<=\langle
\let\>=\rangle

\def\...{...}
\newcommand{\shortStyle}{\textit}
\newcommand{\ie}{\shortStyle{i.e.,}}

\newcommand{\eg}{\shortStyle{e.g.}}

\newcommand{\cf}{\shortStyle{cf.}}

\newcommand{\resp}{resp.}

\makeatletter
\renewcommand*{\eqref}[1]{%
  \hyperref[{#1}]{\textup{\tagform@{\ref*{#1}}}}%
}
\makeatother

\numberwithin{equation}{section}

%% file: sec-homo/introduction.tex
\section{Introduction}
\label{sec:introduction}

%A vertex-transitive (convex) polytope $P\subset\RR^d$ is a polytope in which all vertices are identical under the symmetries of the polytope. More precisely, for any two vertices $v,w\in\F_0(P)$, there is a symmetry $T\in\Aut(G)$, which is an orthogonal transformation of $\RR^d$, that maps $v$ onto $w$, that is, $Tv=w$. Among classes of symmetric polytopes, the class of vertex-transitive polytopes is certainly one of the wildest: Friese and Ladisch \cite{friese2016affine} showed that \emph{almost} any finite group appears as the symmetry group of a vertex-transitive polytope. And finite groups are certainly a wild and not completely understood class of structures.
%
%In this paper we show that vertex-transitivity becomes manageable when restricting to the class of zonotopes.

A \emph{(convex) polytope} is the convex hull of finitely many points. The \mbox{\emph{zonotopes} form} a special class of polytopes for which several equivalent definitions are known (see e.g. \cite[Section 7.3]{ziegler2012lectures}):

\begin{definition}\label{def:zonotope}
A \emph{zonotope} is a polytope $Z\subset\RR^d$ which satisfies any of the follo\-wing equivalent conditions:
\begin{myenumerate}
	\item $Z$ is a parallel projection of a $\delta$-cube.
	\label{it:parallel_projection}
	\item $Z$ is a Minkowski sum of line segments.
	\label{it:minkowski_sum}
	\item $Z$ has only centrally symmetric faces
	\label{it:centrally_symmetric_faces}
	\item $Z$ has only centrally symmetric 2-faces.
	\label{it:centrally_symmetric_2_faces}
\end{myenumerate}
\end{definition}

A polytope $P\subseteq\RR^d$ is called \emph{vertex-transitive} if its (Euclidean) symmetry group $\Aut(P)\subseteq\Ortho(\RR^d)$ acts transitively on its vertex set $\F_0(P)$.
While many classes of symmetric polytopes have been classified in the past (\eg\ regular polyopes, edge- and vertex-transitive polyhedra), a classification of general vertex-transitive polytopes is probably infeasible: 
almost every finite group is isomorphic to the symmetry group of some vertex-transitive polytope
\cite{babai1977symmetry,friese2016affine} (the only exceptions are certain abelian groups and dicyclic groups).
This endeavor can therefore be compared with classifying all finite groups (and their real representations).

%Accepting that general vertex-transitive polytopes can be very wild (\eg\ one thinks about the traveling salesman polytope, representation polytopes, etc.)
We can now impose further restrictions on this class to obtain a new, still interesting, but tractable classification problem.
One well studied sub-class is formed by the \emph{uniform polytopes} (see \eg\ \cite{coxeter1954uniform,johnson1967theory}).
Another large class of polytopes are the \emph{zontopes}, which have seemingly not been probed for their vertex-transitive members before.
%This paper deals with the latter class.
%
In fact, we show that a full classification of vertex-transitive zonotopes can be achieved, and is immediately linked to the classification of finite reflection groups and their root systems.

The first examples of vertex-transitive zonotopes are usually the following: cubes, prisms and permutahedra.
All these belong to the more general class of $\Gamma$-per\-mu\-ta\-he\-dra (see \cref{def:permutahedron}), which contains further examples of vertex-transitive zonopes.
We were able to show, that these are \emph{all} the vertex-transitive zonotopes that exist.
Surprisingly (but in fact easier to prove), the same holds for zonotopes for which all vertices are on a common sphere and all edges are of the same length (no symmetry requirements are necessary).
These will be called \emph{homogeneous} zonotopes.% (see \cref{def:homogeneous}).

As a side product of this classification we obtain further interesting characterizations of root systems.
For example, for a finite centrally symmetric set $R\subset\RR^d\nozero$ consider the intersection of $R$ with a half space that contains exactly half the elements of $R$ (these intersections will be called \emph{semi-stars}, see \cref{def:semi_star_flat}).
For root systems, such \enquote{half-sets} are known as \emph{positive roots}, and it is well-known that the Weyl group acts transitively on these.
We prove an inverse: if in a finite and centrally symmetric set all semi-stars are congruent, then it is a root system.
%We show, that if all these \enquote{half-sets} of $R$ are congruent, then $R$ is a root system.
%For root system, the half-sets are usually called \emph{prositive roots}, and its is well-know that the Weyl group of the roots system acts transitively on them.
%We therefore proved an inverse.
We also prove a further sufficient condition using only the lengths of the sum of the vectors in the semi-stars.

%We will show, that questions and classification attempts concerning regularity and symmetry in zontopes are quite tractable, in contrast to equivalent questions for more general polytopes.
%
%For us, symmetry means \emph{transitivity} on $\delta$-faces, for some fixed $\delta\in\{0,...,d-1\}$, that is, for any two $\delta$-faces $\sigma,\tau\in\F_\delta(Z)$, there is an orthogonal transformation $T\in\Aut(Z)\subset\Ortho(\RR^d)$ that maps $\sigma$ onto $\tau$. 
%We shall focus on $\delta\in\{0,1,2\}$, \ie\ on transitivity on vertices, edges and $2$-faces.
%%Regularity, on the other hand, is a generic name we give to properties that are implied by, \eg\ vertex-transitivity, but can be dealt with using induction arguments.
%Regularity is a more vaguely defined term that generalizes certain symmetry assumptions, and will mean something different from section to section.

%Regularity is defined relative to a symmetry concept, and might look different for each kind of transitivity.
%As a general rule of thumb, if a zonotope is transitive on $\delta$-faces, then it also is regular on these faces. 
%In each case, regularity is constucted to be \emph{hereditary}: if $Z$ is regular on $\delta$-faces, then so are its faces.
%An equivalent statement does not hold for transitivity.
%Regularity is therefore more accessible by induction arguments.

%In this paper, $Z\subset\RR^d$ shall always denote a full-dimensional zonotope. 
%Since zonotopes are always centrally symmetric, we can further assume that $Z$ is centered at the origin.

\subsection{Overview}

%The content of this paper is as follows:
%
\cref{sec:basics} and \cref{sec:permutahedra} recap the relevant definitions and preliminary results that are known from the literature.
In \cref{sec:basics} we define the generators of a zonotope (\ie\ a canonical way to write the zonotope as a Minkowski sum of line segments) and recall how these determine its faces.
We also show that the generators determine the symmetries of the zonotope.
\cref{sec:permutahedra} recalls permutahedra, reflection groups and root systems.
We included proofs for the following facts: permutahedra are exactly the zonotopes generated by root systems.
This is the core property exploited in the following sections.

In \cref{sec:vertex-transitive} we provide a proof of our main result (\cref{res:vertex_transitive_permutahedron}): the vertex-transitive zonotopes are exactly the $\Gamma$-permutahedra, where $\Gamma\subset\Ortho(\RR^d)$ is some finite reflection group.
We also show, that a \emph{homogeneous} zonotope (all vertices on a common sphere, all edges of the same length) is necessarily a $\Gamma$-permutahedron.
These results enable a complete classification of vertex-transitive/homogeneous zonotopes, given in \cref{sec:permutahedron_summary}.

In \cref{sec:root_systems} we apply the classification of the previous sections to give alternative characterizations of root systems (\cref{res:root_halves,res:root_norms}).
%For example, a centrally symmetric set $R\subset\RR^d$ can be intersected with a half space in such a way, that exactly half the vectors of $R$ remain.
%We show, that if all these \enquote{half-sets} of $R$ are congruent, then $R$ is a root system.
%For root system, the half-sets are usually called \emph{prositive roots}, and its is well-know that the Weyl group of the roots system acts transitively on them.
%We therefore proved an inverse.
%We also prove a further sufficient condition.

%% file: sec-homo/basics.tex
\section{Generators, Faces and Symmetries}
\label{sec:basics}

Throughout this paper, $Z\subset\RR^d$ with $d\ge 2$ shall always denote a full-dimensional zonotope.
The cases $d\in\{0,1\}$ are trivial and will not be considered here.
Zonotopes are \emph{centrally symmetric}, and we will assume $Z=-Z$.

\subsection{Generators}

By \cref{def:zonotope} \ref{it:minkowski_sum}, $Z$ is the Minkwoski sum of line segments, and so there is a finite centrally symmetric set $R\subset\RR^d\setminus\{0\}$ with
\begin{equation}\label{eq:sum_of_generators}
Z=\Zon(R):= \sum_{r\in R}\conv\{0,r\} = \Big\{\sum_{r\in R} a_r r\,\Big\vert\, a\in[0,1]^R\,\Big\}.
\end{equation}
We use the convention $\Zon(\eset):=\{0\}$.
%The elements of $R$ are called \emph{generators} of $Z$.
%
%The set of generators of a zonotope is usually not uniquely determined, but there exists a canonical choice:
Conversely,
$$\Gen(Z):=\{r\in\RR^d\mid \text{$\conv\{-r,r\}$ is the translate of an edge of $Z$}\}$$
is a finite, centrally symmetric set with $\Zon(\Gen(Z))=Z$.
The elements of $\Gen(Z)$ will be called \emph{generators} of $Z$.
%
%\begin{myenumerate}
%	\item $0\not\in \Gen(Z)$,
%	\item $\Gen(Z)$ is centrally symmetric, and %\ie\ $\Gen(Z)=-\Gen(Z)$,
%	\item $\Gen(Z)$ is \emph{reduced}, that is, for all $r\in\Gen(Z)$ holds
%	%
%	$$\Span\{r\}\cap \Gen(Z)=\{-r,r\}.$$
%\end{myenumerate}

For this section, let $R\subset\RR^d\setminus \{0\}$ be a finite centrally symmetric set.

\begin{definition}
$R$ is called \emph{reduced} if $\Span\{r\}\cap R=\{-r,r\}$ for all $r\in R$.
%
%Further, a subset $F\subseteq R$ is called a \emph{flat}, if $\Span(F)\cap R=F$.
\end{definition}

The set $\Gen(Z)$ is always \emph{reduced}.
If $R$ is reduced, then it is centrally symmetric, and $R=\Gen(\Zon(R))$, that is, $R$ is the set of generators of some zonotope.

\begin{figure}
\centering
\includegraphics[width=0.6\textwidth]{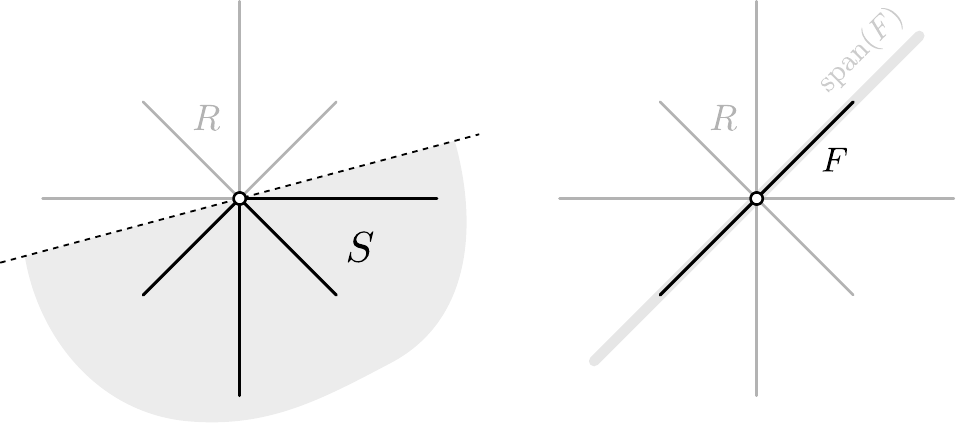}
\caption{Visualization of semi-stars (left) and flats (right).}
\label{fig:flat}
\end{figure}

\begin{definition}\label{def:semi_star_flat}
We consider the following two kinds of subsets of $R$ (see \cref{fig:flat}):
\begin{myenumerate}
	\item A subset $S\subseteq R$ is a \emph{semi-star}, if it contains exactly half the elements of $R$, and can be written as the intersection of $R$ with a half-space.
	\item A subset $F\subseteq R$ is a \emph{flat}, if it can be written as the intersection of $R$ with a linear subspace of $\RR^d$, or alternatively, if $\Span(F)\cap R=F$.
\end{myenumerate}
\end{definition}

The name \enquote{semi-star} is motivated by the occasional use of the term \enquote{star} for centrally symmetric sets $R\subset\RR^d\nozero$.
Note that a semi-star of $R$ contains exactly one element from $\{-r,r\}\subseteq R$ for each $r\in R$.

\subsection{Faces}

It follows a description of how to deduce information about the faces of $Z$ from its generators $\Gen(Z)$.
The proofs of the relevant statements (\cref{res:semi_star_vertices} and \cref{res:linearly_closed_faces}) are found in the standard literature on zonotopes (see e.g. \cite{mcmullen1971zonotopes}, or see the Appendix).

As a general rule of thumb, the shape of the faces is determined by the flats of $\Gen(Z)$, and the position of the vertices is determined by the semi-stars of $\Gen(Z)$.

It follows from \eqref{eq:sum_of_generators} that a vertex $v\in \F_0(Z)$ can be written as 
$$v=\sum_{r\in R} r,\quad\text{for some $R\subseteq\Gen(Z)$}.$$
However, not every subset $R\subseteq\Gen(Z)$ gives rise to a vertex in that way, but only if $R$ is a semi-star of $\Gen(Z)$.
This correspondence is furthermore one-to-one:

\begin{proposition}\label{res:semi_star_vertices}
The vertices $\F_0(Z)$ are in one-to-one correspondence with the semi-stars $S\subset\Gen(Z)$ via the bijection
$S\; \mapsto\;  v_S:=\sum_{r\in S}r.$
\end{proposition}

The faces $\sigma\in\F(Z)$ are themselves zonotopes, and their shapes are determined by their generators.
Their generators are the flats in $\Gen(Z)$:

\begin{proposition}\label{res:linearly_closed_faces}
The following holds:
\begin{myenumerate}
\item 
For each face $\sigma\in\F(Z)$, $\Gen(\sigma)\subseteq\Gen(Z)$ is a flat.
\item
For each flat $F\subseteq\Gen(Z)$ there is a face $\sigma\in\F(Z)$ with $\Gen(\sigma)=F$.
\end{myenumerate}
\end{proposition}

\subsection{Symmetries}

$\Zon(\free)$ and $\Gen(\free)$ are linear, that is, for any orthogonal transformation $T\in\Ortho(\RR^d)$ holds
$$\Zon(TR)=T\Zon(R)\wideand \Gen(TZ)=T\Gen(Z).$$

It immediately follows:

\begin{proposition}\label{res:symmetries}
$\Aut(Z)=\Aut(\Gen(Z))$.
%If $R\subset\RR^d$ is some finite set, then holds $\Aut(R)\subseteq\Aut(\Zon(R))$, with equality if $R$ is reduced.
%
\begin{proof}
If $T\in\Aut(Z)$ is a symmetry of $Z$, then $$T\Gen(Z) = \Gen(TZ)=\Gen(Z),$$
and thus, $T\in\Aut(\Gen(Z))$.
The proof of the other direction goes equivalently.
%
%If $T\in\Aut(R)$ is a symmetry of $R$, then $$T\Zon(R) = \Zon(TR)=\Zon(R),$$
%and thus, $T\in\Aut(\Zon(R))$.
%
%If $R$ is reduced and $0\not\in R$, then $R=\Gen(\Zon(R))$. For all $T\in\Aut(\Zon(R))$
%%
%$$T R= T \Gen(\Zon(R)) = \Gen(T \Zon(R)) = \Gen(\Zon(R)) = R,$$
%%
%and thus, $T\in\Aut(R)$.
%Since $\Aut(R)=\Aut(R\setminus\{0\})$ and $\Zon(R)=\Zon(R\setminus\{0\})$, the result still holds if $0\in R$.
\end{proof}
\end{proposition}

For later use we prove

\begin{proposition}\label{res:semi_star_symmetries}
For $T\in\Ortho(\RR^d)$ holds
\begin{myenumerate}
	\item \label{it:to_semistar}
	if $T\in\Aut(R)$, then $T$ maps semi-stars onto semi-stars.
	\item \label{it:to_symmetry}
	if $TS\subseteq R$ for at least one semi-star $S\subset R$, then $T\in\Aut(R)$.
\end{myenumerate}
\begin{proof}
The semi-stars $S\subset R$ are exactly the subsets of $R$ of the form
$$S=\{r\in R\mid \<r,c\>>0\}$$
for which $c\in\RR^d\nozero$ is non-orthogonal to all $r\in R$.
For $(i)$ compute
\begin{align*}
TS
&=\{Tr\in R\mid \<r,c\>>0\}\\
&=\{r\in TR\mid \<T^{-1}r,c\>>0\}\\%\qquad \qquad \textcolor{lightgray}{\text{using }TR=R}\\
&=\{r\in R\mid \<r,Tc\>>0\} \subset R.
\end{align*}
which is a semi-star since $\<r,Tc\>=\<T^{-1}r,c\>$ cannot be zero, as $T^{-1}r\in R$.

For $(ii)$, assume $TS\subseteq R$.
For any vector $r\in R$, either $r\in S$, and then $Tr\in R$ by assumption, or $-r\in S$, then $Tr=-T(-r)\in -R=R$ by central symmetry.
Thus $TR=R$ and $T\in\Aut(R)$.
\end{proof}
\end{proposition}

%% file: sec-homo/vertex-transitive.tex
\section{Permutahedra and root systems}
\label{sec:permutahedra}

The standard $(d-1)$-dimensional permutahedron is obtained as the convex hull of all coordinate permutations of the vector $(1,2,...,d)\in\RR^{d}$. 
These polytopes, and certain generalizations for general reflection groups (the \emph{$\Gamma$-permutahedra}), provide well-known examples of vertex-transitive zonotopes. 
%The goal of this section is to see that these are already \emph{all} the vertex-transitive zonotopes.
%Equivalently surprising, if $Z$ is a zonotope with all vertices on a common sphere and all edges of the same length (what we shall call a \emph{homogeneous} zonotope), then $Z$ is again already a $\Gamma$-permutahedron.

\subsection{$\Gamma$-permutahedra and reflection groups}

%A first example one usually comes up with when asked about a non-trivial vertex-transitive zonotope, are \emph{permutahedra}. 
%The standard $(d-1)$-dimensional permutahedron is obtained as the convex hull of all coordinate permutations of the vector $(1,2,...,d)\in\RR^{d}$. 
%This is, in fact, a vertex-transitive zonotope. 
%This idea can be generalizes using finite reflection groups.

In the following, $\Gamma\subset\Ortho(\RR^d)$ shall denote a \emph{finite reflection group}, that is, a matrix group %$\Gamma=\<\,T_r\mid r\in R\subset\RR^d\,\>$ generated by reflections 
$$\Gamma=\<\,T_r\mid r\in R\,\>\subset\Ortho(\RR^d),\quad \text{generated by }T_r := \Id-\frac{2rr\T}{\|r\|^2},
%, \qquad (\text{reflection on $r^\bot$})
$$ 
where $R\subset\RR^d\nozero$ is some finite set of vectors, and $T_r\in\Ortho(\RR^d)$ denotes the reflection on the hyperplanes $r^\bot$\!.
Let $\mathcal E$ denotes the union of the reflection hyperplanes $r^\bot,r\in R$.
The connected components of $\RR^d\setminus\mathcal E$ are called \emph{Weyl chambers} of $\Gamma$. 
It is well-known that $\Gamma$ acts regularly (\ie\ transitively and freely) on these chambers \cite{kane2013reflection}.

\begin{definition}\label{def:permutahedron}
%, that is, a group $\Gamma\subseteq\Ortho(\RR^d)$ generated by reflections.
A \emph{$\Gamma$-permutahedron} is a polytope $P\subset\RR^d$ that satisfies any of the following equivalent conditions:
\begin{myenumerate}
	\item $\Gamma$ acts regularly on the vertices of $P$.
	\item $P$ is the \emph{orbit polytope} $\Orb(\Gamma,v):=\conv(\Gamma v)$ of a point $v\in\RR^d\setminus\mathcal E$ in some Weyl chamber of $\Gamma$. %the convex hull of the orbit $\Gamma v$ of $v$ \wrt\ $\Gamma$.
	%, which is a point not fixed by any of the generating reflections $T_r\in \Gamma$.
\end{myenumerate}
% and $v\in \RR^d$ a point not on any of the mirrors of $\Gamma$ (that is, not fixed by any non-trivial $T\in\Gamma$). The convex hull of the orbit $\Gamma v$ is called \emph{$\Gamma$-permutahedron}.\footnote{This construction is also known as the \emph{orbit polytope} $\mathrm{Orb}(\Gamma,v)$.}
\end{definition}

\begin{remark}\label{res:one_vertex_per_Weyl_chamber}
A $\Gamma$-permutahedron has exactly one vertex per Weyl chamber of $\Gamma$, and no vertices in $\mathcal E$ (such a vertex would be fixed by a reflection).
\end{remark}

%A special property of reflection groups is, that a point $v\in\RR^d$ is fixed by some $T\in\Gamma$ if and only if it is already fixed by a generating reflection $T_r\in\Gamma$.

%The similarity to the standard permutahedron lies in the following construction: the $\Gamma$-permutahedron can be obtained as the convex hull of the orbit $\Gamma v$ (a so called \emph{orbit polytope}) when $v\in\RR^d$ is a generic point, that is, a point that is not fixed by any non-trivial $T\in\Gamma$.

%\begin{theorem}\label{res:permutahedra_are_zonotopes}
%$\Gamma$-permuathedra are vertex-transitive zonotopes.
%\end{theorem}

A $\Gamma$-permutahedron is always vertex-transitive, but not necessarily a zonotope (some faces might not be centrally symmetric).
Still, there is always a point $v\in\RR^d\nozero$, so that $\conv(\Gamma v)$ is a zonotope (see, \eg\ the construction in \cref{res:each_group_generates_a_permutahedron}). % (see \cref{res:permutahedron_to_root_system}).
Among these points, there is also always one, so that the zonotope has all edge of the same length.
The resulting polytopes are classically known as the \emph{omnitruncated uniform polytopes}\cite{johnson1967theory}.
%Having all edges of the same length can always be achieved by an appropriate choice of $v$.
%Still another characterization of $\Gamma$-perutahedron is that $\Gamma$ acts transitively and regularly on the vertices of $Z$.

%The standard permutahedron appears among the $\Gamma$-permutahedra as an instance of the $A_d$-permuta\-hedron, with $A_d$ being the reflection group that is also the symmetry group of the regular $d$-simplex (see \cref{fig:permutahedron} for the case $d=3$).

\begin{figure}
\includegraphics[width=0.7\textwidth]{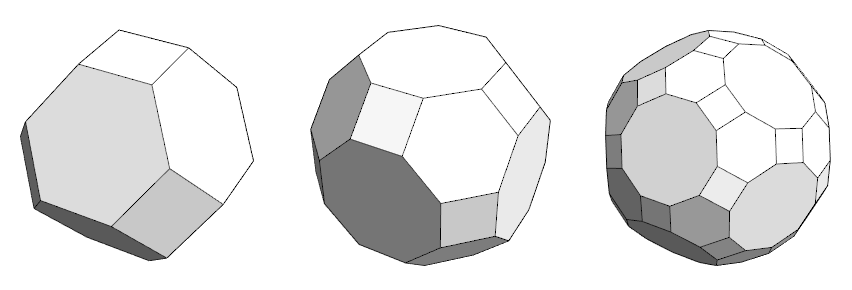}
\caption{$\Gamma$-permutahedron obtained for $\Gamma\in\{A_3, B_3, H_3\}$.}
\label{fig:permutahedron}
\end{figure}

The finite reflection groups are completely classified \cite{kane2013reflection}.
The irreducible ones are denoted by $I_1,I_2(p),A_d,B_d,D_d,H_3,H_4,F_4,E_6,E_7$ and $E_8$ ($d,p\ge 3$, the sub-index always denotes the dimension).
The reducible ones are obtained as direct sums of the irreducible ones.
See \cref{fig:permutahedron} for the $\Gamma$-permutahedra generated from $A_3$ (\resp\ $D_3$), $B_3$ and $H_3$.

%The finite \emph{irreducible} reflection groups are: a single reflection $I_1\cong\Bbb Z_2$ on the line, one parametrized family $I_2(p),p\ge 3$ on the plane (dihedral groups), three groups $A_d$, $B_d$, $D_d$ that exist in all dimension $d \ge 3$, and six exceptional groups $H_3$, $H_4$, $F_4$, $E_6$, $E_7$ and $E_8$ (the subindex always indicates the dimension of the group).
%Additionally, the reducible finite reflection groups are obtained as direct sums of any number of irreducible ones. 

%The goal of this section is to prove that every vertex-transitive zonotopes is a $\Gamma$-permutahedron for some reflection group $\Gamma$.
%The combinatorial type of a $\Gamma$-permutahedron is completely determined by the group $\Gamma$, and as a result of this section, the classification of vertex-transitive zonotopes can be obtained from the classification of finite reflection groups.
%We shall discuss this correspondence in detail in \cref{sec:permutahedron_summary}.
%This result is in stark contrast to the case of general vertex-transitive polytopes, for which a classification is infeasible: almost all finite groups are symmetry groups of vertex-transitive polytopes (see \cite{babai1977symmetry,friese2016affine}).

\subsection{Root systems}

Root systems bridge between zonotopes and reflection groups.
We shall work with a minimalistic version of the definition (no crystallographic or reducedness restrictions):

\begin{definition}
%A finite set $R\subset\RR^d$ is called a \textit{root system}, if it is invariant \wrt\ the reflection $T_r$ for all non-zero $r\in R$.
A \emph{root system} is a finite set $R\subset\RR^d\nozero$, for which $T_r R=R$ for all vectors $r\in R$.
%A \textit{root system} is a finite set of vectors $\R\subset\RR^d$ with the following properties:
%%
%\begin{myenumerate}
%	\item for every $r\in\R$ we have $\R\cap \Span\{r\}=\{-r,r\}$, and
%	\item for every $r\in\R$, $\R$ is invariant under reflection on $r^\bot$ (the hyperplane through the origin, and orthogonal to $r$).
%\end{myenumerate}
\end{definition}

We can give an alternative definition of $\Gamma$-permutahedron using root systems: a zonotope is a $\Gamma$-permutahedron if and only if $\Gen(Z)$ is a root system (we provide a proof for completeness in \cref{res:root_system_to_permutahedron} and \cref{res:permutahedron_to_root_system}).

The \emph{Weyl group} of a root system $R$ is the group 
$$\Gamma :=\<\,T_r\mid r\in R\,\>.$$
This is a finite reflection group. 
Note that $\Gamma\subseteq\Aut(R)$ by definition of root system.

\begin{lemma}\label{res:root_system_to_permutahedron}
%If $R\subset\RR^d$ is a root system, then $\Zon(R)$ is a $\Gamma$-permutahedron with reflection group $\Gamma:=\<\,T_r\mid r\in R\,\>$.

If $\Gen(Z)$ is a root system, then $Z$ is a $\Gamma$-per\-mu\-ta\-he\-dron, where $\Gamma$ is the Weyl group of $\Gen(Z)$.
\begin{proof}
%Assume that $\Star(Z)$ is a root system.
%Define the finite reflection group 
%%
%$$\Gamma:=\<\,T_r\mid r\in\Star(Z)\,\>.$$
%%

It holds $\Gamma\subseteq\Aut(\Gen(Z))$ by definition of root system, and therefore $\Gamma\subseteq\Aut(Z)$ by \cref{res:symmetries}.

The proof proceeds in two steps: first, we show that $\Gamma$ acts transitively on $\F_0(Z)$.
This shows that $\F_0(Z)$ is an orbit of $\Gamma$.
Second, we show that no $T_r,r\in\Gen(Z)$ fixes a vertex of $Z$.
The vertices must then lie in the interior of the Weyl chambers of $\Gamma$, and $\Gamma$ must act regularly on them.

%, and second, we show that $\Gamma$ acts freely, by showing that no $T_r,r\in\Gen(Z)$ fixes a vertex.

%We first show that $Z$ is vertex-transitive.
\emph{Step 1}: The edge graph of $Z$ is connected. So, for any two vertices $v,w\in \F_0(Z)$ the edge graph contains a path $v_0v_1\cdots v_k$ between $v=v_0$ and $w=v_k$.
That means, $e_i:=\conv\{v_{i-1},v_i\}$ is an edge of $Z$ for all $i\in\{1,\...,k\}$.
Choose $r_i\in \Gen(Z)$ parallel to $e_i$.
The reflection $T_{r_i}$ fixes $e_i$ and exchanges its end vertices $v_{i-1}$ and $v_i$.
The map $T:=T_{r_k}\!\cdots T_{r_1}\in\Gamma$ now satisfies $Tv=w$.
This proves vertex-transitivity.

%%Then $Z$ and $\Star(Z)$ are invariant \wrt\ $\Gamma$.
%Let $e\in\F_1(Z)$ be an edge and $r\in\Star(Z)$ so that $e=\conv\{-r,r\}$. 
%Thus $T_r$ fixes $e$ set-wise and exchanges its end vertices.
%Since the edge-graph of $Z$ is connected, for any two vertices $v,w\in\F_0(Z)$ there is a path $v=v_0v_1\cdots v_k=w$ in the edge graph.
%If the edge $v_{i-1}v_{i}=\conv\{-r_i,r_i\}$ for $r_i\in\Star(Z)$, then $T:=T_{r_1}\cdots T_{r_{k}}\in\Gamma$ maps $v$ onto $w$, and we obtained vertex-transitivity.
%
%%a sequence of such reflections$\Gamma$ acts transitively on the vertices of $Z$.

%We prove that the $T_r,r\in\Gen(Z)$ act freely on the vertices of $Z$.
\emph{Step 2}: 
%Suppose that $T_r$ fixes $v\in\F_0(Z)$.
%By \cref{res:semi_star_vertices}, there exist a \emph{unique} semi-star $\Gen_c(Z)$ whose elements sum up to $v$.
%W.l.o.g.\ assume $r\in\Gen_c(Z)$.
%Since $T_r r=-r$ we have $r\not\in T_r\Gen_c(Z)$. %, and the reflection $T_r$ does not fix $\Gen_c(Z)$ set-wise.
%This means $T_r\Gen_c(Z)$ is another semi-star that sums up to $v$, in contradiction to its uniqueness.
%
For $v\in\F_0(Z)$, there is a \emph{unique} semi-star $S\subset \Gen(Z)$ whose elements sum up to $v$ (\cref{res:semi_star_vertices}).
For any generator $r\in \Gen(Z)$ \emph{exactly} one element of $\{-r,r\}$ is in $S$.
Since $T_r (\pm r)=\mp r$, the reflection $T_r$ cannot fix $S$ but maps it to a different semi-star (by \cref{res:semi_star_symmetries} \ref{it:to_semistar}). 
Since the semi-star is unique for $v$, $T_r$ cannot fix $v$ either.
\end{proof}
\end{lemma}

\begin{remark}\label{res:each_group_generates_a_permutahedron}
For any finite reflection group $\Gamma$, there exists a $\Gamma$-permutahedron which is a (vertex-transitive) zonotope.

To see this, choose a \emph{reduced} root system $R$ that generates $\Gamma$ (\eg\ build $R$ from the normals of the mirror hyperplanes of $\Gamma$).
Then, $\Zon(R)$ is a zonotope, and a $\Gamma$-permutahedron by $\Gen(\Zon(R))=R$ and \cref{res:root_system_to_permutahedron}.
\end{remark}

%\begin{lemma}\label{res:permutahedron_to_root_system}
%%If $Z$ is a $\Gamma$-permutahedron, then $\Gen(Z)$ is a root system.
%If $R\subset \RR^d$ is a finite and centrally symmetric set, so that $Z:=\Zon(R)$ is a $\Gamma$-permutahedron, then the reduction of $R$ is a root system.
%%
%\end{lemma}
%
%Recall the reduction of $R$, which is the reduced set $\Gen(\Zon(R))$.
%The set $R$ itself does not have to be a root system.
%Setting $R:=\Gen(Z)$ for a $\Gamma$-permutahedron $Z$, would show that $\Gen(Z)$ is a root system.

\begin{lemma}\label{res:permutahedron_to_root_system}
If $Z$ is a $\Gamma$-permutahedron, then $\Gen(Z)$ is a root system.
%Let $R\subset \RR^d$ be a finite and centrally symmetric set, and $Z:=\Zon(R)$.
%If $\Aut(Z)\subseteq\Aut(R)$ and $Z$ is a $\Gamma$-permutahedron, then $R$ is a root system.
%
\begin{proof}%[Proof of \cref{res:permutahedron_to_root_system}]
%Choose a vector $r\in R$ and let $e\in\F_1(Z)$ be an edge parallel to it.
%% (which exists by \cref{res:linearly_closed_faces}, \msays{not exactly}).
%Since the end vertices of $e$ cannot be in the same Weyl chamber of $\Gamma$, $e$ must cross one of the reflection hyperplanes of $\Gamma$.
%These are symmetry hyperplanes of $Z$, and since $e$ crosses one, it must be perpendicular to it.
%Thus, this hyperplane is $r^\bot$\!.
Choose a generator $r\in \Gen(Z)$ and let $e\in\F_1(Z)$ be an $r$-parallel edge.
% (which exists by \cref{res:linearly_closed_faces}, \msays{not exactly}).
There is at most one vertex per Weyl chamber  of $\Gamma$ (\cref{res:one_vertex_per_Weyl_chamber}), hence the end vertices of $e$ cannot be in the same chamber.
Therefore, $e$ must cross one of the reflection hyperplanes of $\Gamma$.
These are symmetry hyperplanes of $Z$, and if $e$ crosses one, it must be perpendicular to it.
Thus, this hyperplane is $r^\bot$\!.

We have shown that $T_r\in\Aut(Z)=\Aut(\Gen(Z))$ for all $r\in R$, and therefore, $\Gen(Z)$ is a root system.
\end{proof}
\end{lemma}

%The formulation of \cref{res:permutahedron_to_root_system} is more general than necessary for this section.
%For immediate use, set $R:=\Gen(Z)$, and \cref{res:symmetries} ensures that $\Aut(\Gen(Z))=\Aut(Z)$.
%This more general version will be used in \cref{sec:edge-transitive}.

%\begin{corollary}\label{res:permutahedron_to_root_system_cor}
%If $Z$ is a $\Gamma$-permutahedron, then $\Gen(Z)$ is a root system.
%\end{corollary}

%The cases of $d$-dimensional vertex-transitive zontopes with $d=0$ and $d=1$ are trivial: the only polytope of these dimensions are the single vertex and an interval, both are vertex-transitive zonotopes.
%In the following, we restrict to $d\ge 2$.

\section{Vertex-transitive and homogeneous zonotopes}
\label{sec:vertex-transitive}

The goal of this section is to prove \cref{res:vertex_transitive_permutahedron}, that all vertex-transitive (and homogeneous, \cref{def:homogeneous}) zonotopes are $\Gamma$-permutahedra.
%From this, we derive a full classification in \cref{sec:permutahedron_summary}.

\subsection{Roadmap}
We briefly describe the roadmap to the proof.

The idea is to prove the statement for 2-dimensional zonotopes (\cref{sec:2d}), and then, transfer this result to general dimensions (\cref{sec:general_case}).
This is possible by the following result (\cref{res:2face_root_system}): if all 2-faces of a zonotope are vertex-transitive, then it is a $\Gamma$-permutahedron.

%\cref{sec:2d} discusses 2-dimensional zonotopes, \ie\ we argue that 2-dimen\-sional vertex-transi\-tive zonotopes are $\Gamma$-permutahedra (where $\Gamma\in\{I_2(p),I_1\oplus I_2\}$ is a 2-dimensional reflection group).
%This case is quite trivial, but it serves as a basis for the general result: we can now prove that a zonotope is already a $\Gamma$-permutahedron, if all its 2-faces~are vertex-transitive (\cref{res:2face_root_system}).

We cannot immediately see, that the 2-faces of our vertex-transitive zonotopes are vertex-transitive (this is true, but not obvious, and it is false for general polytopes).
We therefore proof a helper result (\cref{res:2d_vertex_transitive}) which further weakens the condition of vertex-transitivity of the 2-faces.

%The faces of a vertex-transitive polytope are \emph{not} necessarily vertex-transitive.
%We therefore cannot apply \cref{res:2face_root_system} directly to our vertex-transitive zonopes (their faces are vertex-transitive, but we cannot prove it at this point).
%Instead, we establish a helper result (\cref{res:2d_vertex_transitive}) which further weakens the condition of vertex-transitivity of the 2-faces.

In \cref{sec:general_case} we define \emph{homogeneous} zonotopes (all vertices on a sphere, all edges of the same length), and prove that these are $\Gamma$-per\-mu\-ta\-he\-dra (\cref{res:homogeneous_root_system}).
The proof is surprisingly simple: their 2-faces are regular polygons, hence vertex-transitive.
To transfer the result to vertex-transitive zonotopes, we define the \emph{normalization} of a zonotope (\cref{def:normalization}), which transforms any vertex-transitive zonotope into a homogeneous one.
This transformation preserves the edge-directions of the zonotope, so that with the help of \cref{res:2d_vertex_transitive} we can show that all 2-faces of the initial zonotope must have been vertex-transitive.
%Therefore, \cref{res:2d_vertex_transitive} can be used to show, that the 2-faces of the initial vertex-transitive zonotope must have been vertex-transitive.
Applying \cref{res:2face_root_system} then proves the main result \cref{res:vertex_transitive_permutahedron}.

\subsection{2-dimensional zonotopes}
\label{sec:2d}

%The roadmap for this section is simple: show that $\R(Z)$ for a vertex-transitive zonotope is a root system, and then show that a zonotope with such roots is a $\Gamma$-permutahedron.

\begin{figure}
\centering
\includegraphics[width=0.65\textwidth]{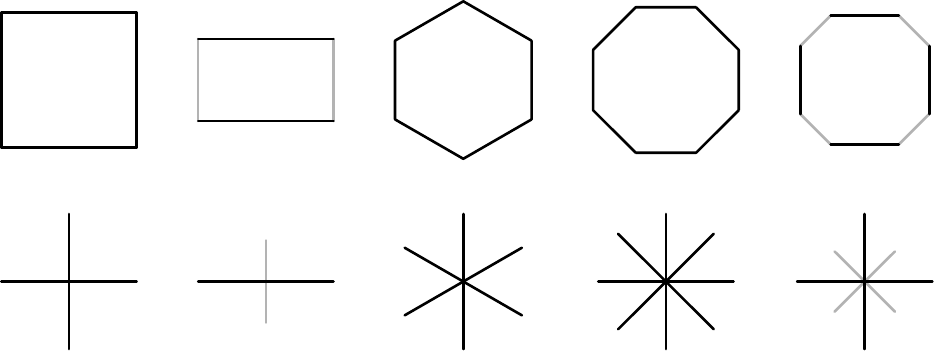}
\caption{Some $2$-dimensional vertex-transitive zonotopes and their generators (not necessarily at scale).}
\label{fig:2d}
\end{figure}

%We first investigate 2-dimensional vertex-transi\-ti\-ve zonotopes as an important special case. 
A 2-dimensional zonotope is a centrally symmetric $2n$-gon. 
Such one is vertex-transitive, if either
\begin{myenumerate}
	\item it is a regular $2n$-gon, or
	\item $n$ is even, and it has alternating edge lengths, as seen in \cref{fig:2d}.
\end{myenumerate}
This list is complete: every vertex-transitive polygons is an orbit polytope to a~cyclic group or dihedral group. These orbit polytopes are contained in the list above.

If $Z$ is a $2n$-gon as listed in $(i)$ or $(ii)$, $\Gen(Z)$ consists of of $2n$ vectors in $\RR^2$, equally spaced by an angle of $\pi/n$, and in the case $(ii)$, these vectors alternate in length (see \cref{fig:2d}).
These are exactly the root systems that corresponds to the reflection group $I_2(n)$ (if $n\ge 3$), or  $I_1\oplus  I_1$ (if $n=2$).
%These are all possible 2-dimensional reflection groups.
%The root system consists of $2n$ vectors in $\RR^2$, equally spaces by an angle of $\pi/n$, and if $n$ is even, of possibly alternating lengths (\cf\ \cref{fig:2d}).
%Consequently

Applying \cref{res:root_system_to_permutahedron}, we obtain the classification in dimension two.

\begin{corollary}\label{res:dim_2}
A $2$-dimensional vertex-transitive zonotope is a $\Gamma$-per\-mu\-ta\-he\-dron.
\end{corollary}

This already finishes the case of 2-dimensional vertex-transitive zonotopes.
This case is important for the following reason: investigating vertex-transitive zonotopes in general dimensions comes down to the study of their 2-faces:

\begin{theorem}\label{res:2face_root_system}
If all $2$-faces of $Z$ are vertex-transitive, then $Z$ is a $\Gamma$-permutahedron.
\begin{proof}
Choose generators $r,s\in\Gen(Z)$.
We show $T_{r}s\in\Gen(Z)$, establishing that $\Gen(Z)$ is a root system and $Z$ is a $\Gamma$-permutahedron by \cref{res:root_system_to_permutahedron}.
%This is trivial if $r_1$ and $r_2$ are linearly dependent.

The case $r=\pm s$ is trivial.
We therefore assume that $$R:=\Gen(Z)\cap\Span\{r,s\}$$ is 2-dimensional.
%We therefore assume that $r$ and $s$ are linearly independent.
In particular, $R\subseteq\Gen(Z)$ is a 2-dimensional flat of $\Gen(Z)$.
By \Cref{res:linearly_closed_faces} there exist a 2-face $\sigma\in\F_2(Z)$ with $\Gen(\sigma)=R$.
By assumption, $\sigma$ is vertex-transitive, and $R$ therefore a root system (\cref{res:dim_2}).
Conclusively, $T_r s\in R\subseteq \Gen(Z)$.
\end{proof} 
\end{theorem}

In order to apply \cref{res:2face_root_system}, we need to prove that certain 2-faces are vertex-transitive. 
This does not follows immediately (even though it is true) from the fact that we start from a vertex-transitive zonotope.
Instead, we use the following helper result:

\begin{proposition}\label{res:2d_vertex_transitive}
%Let $Z$ be a 2-dimensional zonotope which
If $Z\subset\RR^2$ is a $2$-dimensional zonotope which
\begin{myenumerate}
	\item has all vertices on a common circle, and
	\item has the same edge directions as a regular $2n$-gon,
	% (but not necessarily the same edge lengths) as a regular $2n$-gon,
\end{myenumerate}
then $Z$ is vertex-transitive.
\end{proposition}

This statement is elementary.
We sketch its proof:

\begin{figure}
\centering
\includegraphics[width=0.3\textwidth]{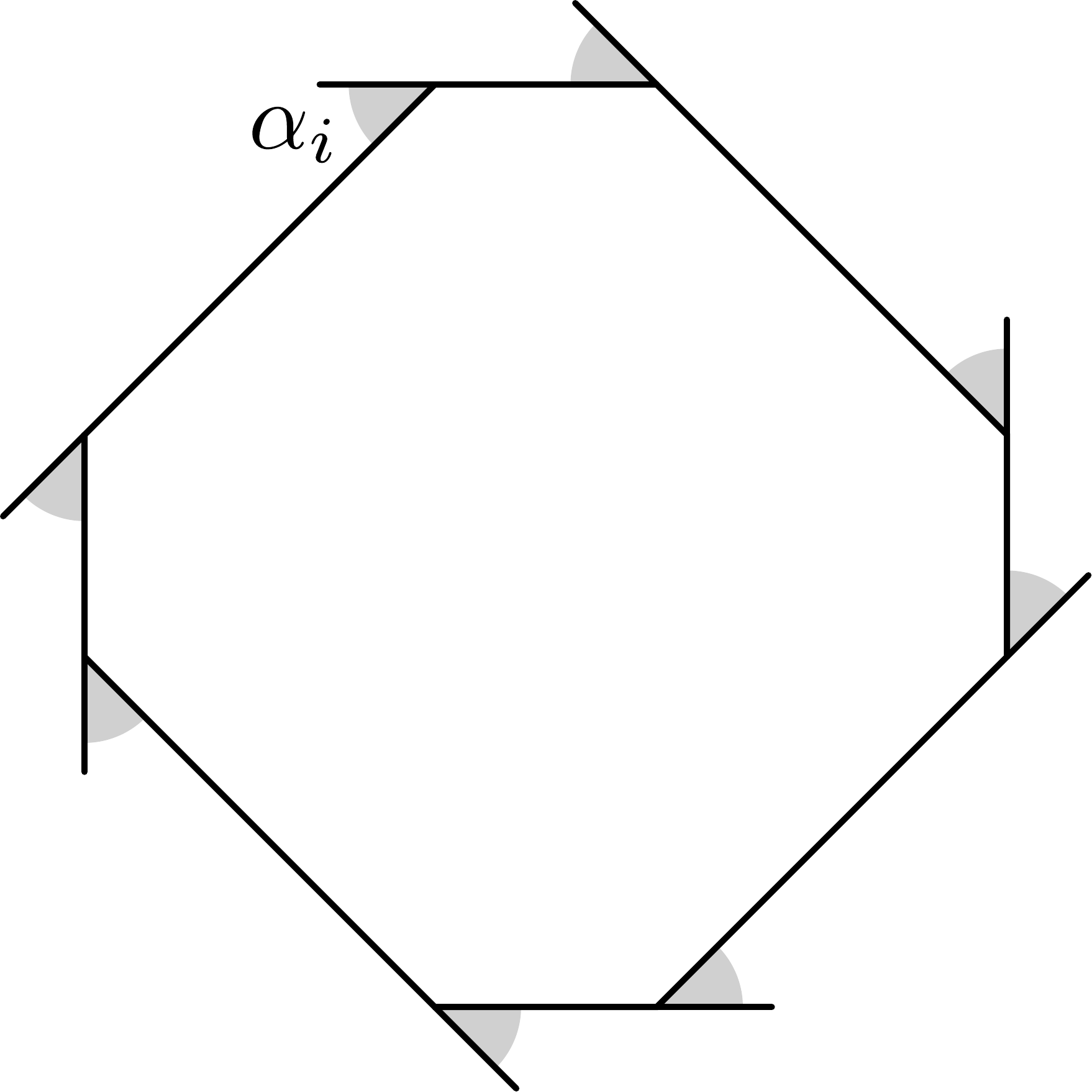}
\caption{The exterior angles of a (convex) polygon always add up to $2\pi$ (used in the proof of \cref{res:2d_vertex_transitive}).}
\label{fig:exterior_angles}
\end{figure}

\begin{proof}[Proof of \cref{res:2d_vertex_transitive}]
A (convex) polygon has at most two edges of the same direction, and if it is centrally symmetric, it has \emph{exactly} two of each.
A regular $2n$-gon has $n$ edge directions, and by $(ii)$, so has $Z$.
Since $Z$ is centrally symmetric, it must be a $2n$-gon as well.
%Hence, if $Z$ is a $2m$-gon, it has exactly $m$ edge directions. 
%By $(ii)$, these directions match the ones of a regular $2n$-gon, which has $n$ directions, and so $m=n$ and $Z$ is a $2n$-gon.

Let $\alpha_i\in\RR$ be the exterior angle of the $i$-th vertex of $Z$ (see \cref{fig:exterior_angles}).
By $(ii)$ we have $\alpha_i = k_i \pi/n$, where $k_i\in\mathbb N$ is an integer $\ge 1$.
The exterior angles of a (convex) polygon add up to $2\pi$, and so we estimate
$$2\pi = \sum_{i=1}^{2n} \alpha_i = \sum_{i=1}^{2n} \frac{k_i \pi}n \ge 2n\cdot \frac{1\cdot \pi}n=2\pi,$$
and conclude that $k_i=1$. 
That is, all exterior and \emph{interior} angles are equal.

For simplicity, assume that $Z$, and all poygons mentioned in the following paragraph are of circumradius one.
Let $\ell$ be the length of the shortest edge of $Z$.
Then, $\ell$ is no longer than the edge length of a regular $2n$-gon.
This ensures that there exists a vertex-transitive $2n$-gon with an edge of length $\ell$ (consider an appropriately chosen orbit polytope of $I_2(n)$ \resp\ $I_1\oplus I_1$).
But a polygon satisfying $(i)$ and with prescribed identical interior angles at every vertex is already uniquely determined by placing a single edge (the placement of both incident edges follows uniquely from the set restrictions, and this iteratively determines the whole polygon).
Therefore, $Z$ must be this vertex-transitive polygon. %\TODO
\end{proof}

%Now to the general case, which is not much harder, as we just reduce to the 2-dimensional case. However, we need the following lemma:
%
%\begin{lemma}
%If $Z$ is a vertex-transitive zonotope and $\bar Z$ the zonotope generated by $$\bar\R:=\left\{\frac r{\|r\|} \;\middle\vert\; r\in\R(Z)\right\},$$ then $\bar Z$ is also vertex-transitive.
%%
%\begin{proof}
%...
%\end{proof} 
%\end{lemma}

\subsection{The general case}
\label{sec:general_case}

%We now introduce \emph{homogeneous} zonotopes and show that these must be $\Gamma$-permutahedra.
%We then use this result to proof the main theorem \cref{res:vertex_transitive_permutahedron}.

\begin{definition}\label{def:homogeneous}
A zonotope is said to be \emph{homogeneous}, if all its vertices are on a common sphere, and all its edges are of the same length.
\end{definition}

Homogeneity is a \emph{hereditary} property:

\begin{observation}\label{res:homogeneous}
The faces of a homogeneous zonotope $Z$ are homogeneous: all edges of a face $\sigma\in\F(Z)$ are of the same length.
All vertices of $Z$ are on a sphere $\mathcal S$, and all vertices of $\sigma$ are on an affine subspace of $\RR^d$.
All vertices of $\sigma$ are on the intersection of this subspace with $\mathcal S$, which is itself a sphere.
\end{observation}

The 2-faces of homogeneous zonotopes are homogeneous, and homogeneous 2-faces are regular, thus vertex-transitive. With \cref{res:2face_root_system} we conclude

\begin{corollary}\label{res:homogeneous_root_system}
If $Z$ is homogeneous, then $Z$ is a $\Gamma$-permutahedron.
\end{corollary}

To apply this to vertex-transitive zonotopes via the following construction:

\begin{definition}\label{def:normalization}
The \emph{normalization} of $Z$ is the zonotope 
$$Z^*:=\Zon\Big(\Big\{\frac r{\|r\|}\;\Big\vert\; r\in\Gen(Z)\Big\}\Big).$$
\end{definition}

The normalized zonotope has the same edge directions as $Z$, but all edges are of the same length.
%
%Each face $\sigma\in\F(Z)$ is congruent to $\Zon(R)$ for some linearly closed subset $R\subseteq\Star(Z)$. The corresponding face $\sigma^*$

%\begin{definition}
%Let $R\subseteq\RR^d$ be a set and $c\in\RR^d$ a vector that is not perpendicular to any vector of $R$.
%The set
%%
%$$R_c:=\{r\in R\mid\<r,c\>>0\}$$
%%
%is called a \emph{half} of $R$.
%\end{definition}
%
%\begin{observation}
%Each half is a sign selection of $\Star(Z)$.
%This gives us a one-to-one correspondence between the halves $\Star_c(Z)$ and the vertices of $Z$, namely,
%%
%$$\Star_c(Z) \;\mapsto\; v_c:=\sum_{\mathclap{r\,\in\, \Star_c(Z)}} r.$$
%%
%Any vector $c$ that defines a half $\Star_c(Z)$, is the normal vector of the hyperplane that touches $Z$ in the vertex $v_c$.
%\end{observation}

We further need to understand how vertex-transitivity is determined by $\Gen(Z)$: 

\begin{proposition}\label{res:half_properties}
$Z$ is vertex-transitive if and only if all semi-stars of $\Gen(Z)$~are congruent.
\begin{proof}
Consider two vertices $v,v'\in\F_0(Z)$ determined by semi-stars $S,S'\subseteq \Gen(Z)$ (via \cref{res:semi_star_vertices}).

A symmetry $T\in\Aut(Z)$ that maps $v$ onto $v'$, maps $S$ onto a \mbox{semi-star} whose elements add up to $v'$ (\cref{res:semi_star_symmetries}).
But the semi-star with this property is unique by \cref{res:semi_star_vertices}, and therefore must equal $S'$.
Thus, $TS=S'$, and $T$ is a congruence between the semi-stars.

On the other hand, if $T\in\Ortho(\RR^d)$ maps the semi-star $S$ onto $S'\subset\Gen(Z)$, then it also maps $v$ onto $v'$. % (\cref{res:symmetries}).
And $T\in\Aut(Z)$ by \cref{res:semi_star_symmetries} $(ii)$.
\end{proof}
\end{proposition}

%As an easy but relevant consequence, we have

\begin{corollary}\label{res:norm_transitive}
If $Z$ is vertex-transitive, then so is its normalization $Z^*$.
\begin{proof}
If all semi-stars of $\Gen(Z)$ are congruent, then this certainly stays valid when normalizing all vectors in $\Gen(Z)$.
%We show that all semi-stars of $\Gen(Z^*)$ are congruent.
%Vertex-transitivity then follows with \cref{res:half_properties}.
%
%Choose $c,c'\in\RR^d$ in general position \wrt\ $\Gen(Z)$.
%%Each half of the normalization, \eg\ $\Star_c(Z^*)$, can be written as
%Note that
%%
%$$\Gen_c(Z^*) = \Big\{\frac r{\|r\|} \,\Big\vert\, r\in\Gen_c(Z) \Big\}.$$
%%
%and equivalently for $c'$.
%But $\Star_c(Z)$ is congruent to $\Star_{c'}(Z)$ via some transformation $T\in\Aut(Z)$.
%We then find
%%
%\begin{align*}
%T\Star_c(Z^*) 
%&= \Big\{\frac {Tr}{\|r\|} \,\Big\vert\, r\in\Star_c(Z) \Big\} = \Big\{\frac {Tr}{\|Tr\|} \,\Big\vert\, r\in\Star_c(Z) \Big\}\\
%&= \Big\{\frac {r}{\|r\|} \,\Big\vert\, r\in T\Star_c(Z) \Big\}= \Big\{\frac {r}{\|r\|} \,\Big\vert\, r\in \Star_{c'}(Z) \Big\}= \Star_{c'}(Z^*).
%\end{align*}
\end{proof}
\end{corollary}

%We mention the following immediate consequence:

\begin{corollary}\label{res:norm_homogeneous}
If $Z$ is vertex-transitive, then $Z^*$ is a $\Gamma$-per\-mu\-ta\-he\-dron.
\begin{proof}
By \cref{res:norm_transitive}, $Z^*$ is vertex-transitive, in particular, all vertices are on a common sphere.
Furthermore, all edges of $Z^*$ have the same length.
Thus, $Z^*$ is homogeneous, and a $\Gamma$-permutahedron by \cref{res:homogeneous_root_system}.
\end{proof}
\end{corollary}

We can finally prove the main result:

\begin{theorem}\label{res:vertex_transitive_permutahedron}
A vertex-transitive zonotope $Z$ is a $\Gamma$-permutahedron.
\begin{proof}
In \cref{res:norm_homogeneous} we saw that $Z^*$ is a $\Gamma$-permu\-ta\-he\-dron, and $\Gen(Z^*)$ therefore a root system by \cref{res:root_system_to_permutahedron}.
For each 2-face $\sigma\in\F_2(Z)$, the set
$$R^*:=\Big\{\frac r{\|r\|}\,\Big\vert\,r\in\Gen(\sigma) \Big\} = \Gen(Z^*)\cap \Span(\Gen(\sigma))$$
is a 2-dimensional flat subset of the root system $\Gen(Z^*)$.
As such, it is a root system itself.
A 2-dimensional root system consists of vectors that are equally spaces by an angle $\pi/n$ for some $n\in \NN$, or in other words, the elements of $R^*$ are the edge directions of a regular $2n$-gon.
The 2-face $\sigma$ has the edge direction (but not necessarily the edge lengths) from $R^*$, hence, the same edge direction as a regular $2n$-gon.
%Hence, $\sigma^*:=\Zon(R^*)$ is a $\Gamma$-permutahedron.
%As $\sigma$ has the same edge directions as the $\Gamma$-permutahedron $\sigma^*$, \cref{res:2d_vertex_transitive} states that $\sigma$ is vertex-transitive.
Additionally, as a face of a vertex-transitive polytope, $\sigma$ has all vertices on a common sphere. By \cref{res:2d_vertex_transitive}, $\sigma$ is vertex-transitive.

We found that all 2-faces of $Z$ are vertex-transitive. 
\Cref{res:2face_root_system} provides that $Z$ is a $\Gamma$-per\-mu\-tahedron.
\end{proof}
\end{theorem}

We summarize the results in the following theorem:

\begin{theorem}\label{res:zonotope_equivalents}
If $Z$ is a zonotope, then the following are equivalent:
\begin{myenumerate}
	\item $Z$ is vertex-transitive.
	\item all semi-stars of $\Gen(Z)$ are congruent.
	\item $\Gen(Z)$ is a root system.
	\item $Z$ is a $\Gamma$-permutahedron.
\end{myenumerate}
\begin{proof}
%We have previously shown 
%$(ii)\Leftrightarrow(i)\Leftrightarrow(iv)\Leftrightarrow(iii)$.
$
(ii)
\;\smash{\overset{\raisebox{0.3ex}{\tiny\ref{res:half_properties}}}\Longleftrightarrow}\;
(i)
\;\smash{\overset{\raisebox{0.3ex}{\tiny\ref{def:permutahedron}+\ref{res:vertex_transitive_permutahedron}}}\Longleftrightarrow}\;
(iv)
\;\smash{\overset{\raisebox{0.3ex}{\tiny\ref{res:root_system_to_permutahedron}+\ref{res:permutahedron_to_root_system}}}\Longleftrightarrow}\;
(iii).
$
\end{proof}
\end{theorem}

We list some consequences of this classification, none of which holds for general polytopes:

\begin{corollary}
The following holds:
\begin{myenumerate}
	\item The faces of a vertex-transitive zonotope are vertex-transitive.
	\item A homogeneous zonotope (\ie\ it has all vertices on a common sphere, and all edges of the same length) is vertex-transitive.
	%\item A vertex-transitive zonotope, which is not combinatorially equivalent to the $d$-cube, has a homogeneous facet.
	\item A zonotope in which all faces (in fact, all $2$-faces) are homogeneous (\resp\ vertex-transitive) is itself homogeneous (\resp\ vertex-transitive).
\end{myenumerate}
\begin{proof}
If $Z$ is vertex-transitive, then $\Gen(Z)$ is a root system (\cref{res:vertex_transitive_permutahedron}). By \cref{res:linearly_closed_faces}, for every face $\sigma\in\F(Z)$, its generators $\Gen(\sigma)\subseteq\Gen(Z)$ are a flat subset of the root system $\Gen(Z)$, hence, a root system as well.
Hence, $\sigma$ is a $\Gamma$-permutahedron and vertex-transitive.
This proves $(i)$.

Part $(ii)$ follows immediately from \cref{res:homogeneous_root_system}. Part $(iii)$ follows from \cref{res:2face_root_system} (a homogeneous 2-face is regular, hence also vertex-transitive).
\end{proof}
\end{corollary}

%Counter-examples for general polytopes are (in this order): certain non-uniform anti-prisms (see \cref{fig:...}), elongated square gyrobicupola (which looks almost like the rhombicuboctahedron), the anti-prism again, and a square-based pyramid (in fact, most Johnson solids).

\subsection{The classification of vertex-transitive/homogeneous zonotopes}
\label{sec:permutahedron_summary}

We apply the results in \cref{res:homogeneous_root_system} and \cref{res:vertex_transitive_permutahedron} to compile a list of all vertex-transitive/homogeneous zonotopes.
%Surely, this list is easily obtained by using the classification of finite reflection groups and their root systems.

%We explicitly describe the vertex-transi\-ti\-ve zonotopes by deducing them from the classification of finite reflection groups and their root systems.

%We restrict to the vertex-transitive zonotopes in $d \ge 2$ dimensions.
%For $d=0$ we have the single vertex, and for $d=1$ we have the interval.

We restrict the enumeration to the \emph{irreducible} cases, that is, those which result from irreducible reflection groups.
The remaining zonotopes (the reducible ones) are obtained as cartesian products of the irreducible ones.
The most prominent members of this class are the hypercubes and prisms.

We obtain the following list of irreducible \emph{homogeneous} zonotopes:

\begin{myenumerate}
\item 
infinitely many 2-dimensional homogeneous zonotopes (the regular $2n$-gons),
\item
for each $d\ge 3$, the $d$-dimensional zonotopes generated by reflection groups $A_d, B_d$ and $D_d$ (which are distinct if and only if $d>3$), and
\item
six exceptional zonotopes to the reflection groups $H_3$, $H_4$, $F_4$, $E_6$, $E_7$ and $E_8$ in their respective dimensions $d\in\{3,4,6,7,8\}$.
\end{myenumerate}

All of these are uniquely determined up to scale and orientation.
The polytopes in that list are also classically known as the \emph{omnitruncated uniform polytopes}.
This terminology, and many related names were coined by Norman Johnson \cite{johnson1967theory}.

\begin{figure}
\centering
\includegraphics[width=0.9\textwidth]{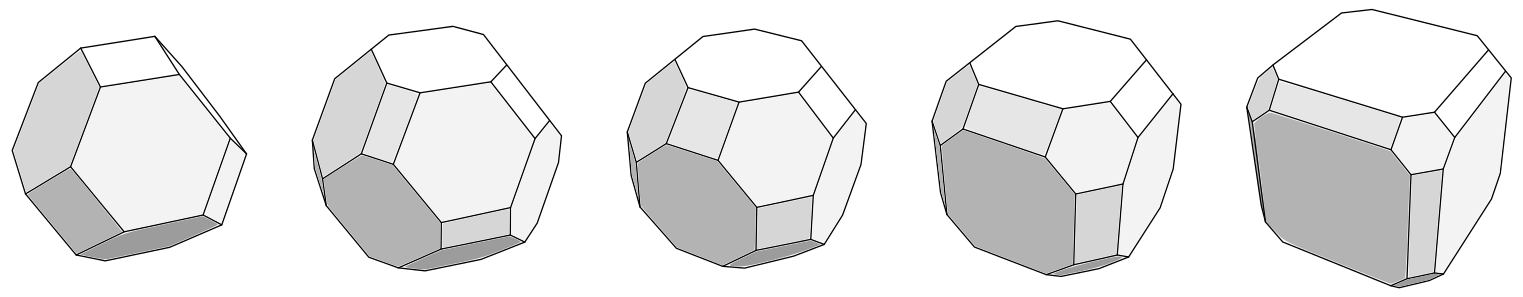}
\caption{Sample from the 1-dimensional family of $B_3$-zonotopes, except for the left one, which degenerated to the $A_3/D_3$-zonotope. The middle one is the homogeneous representative (\emph{omnitruncated cube}).}
\label{fig:B3_mod}
\end{figure}

The (irreducible) \emph{vertex-transitive} zonotopes differ from the homogeneous ones in three cases:
the zonotopes that correspond to the root systems $I_2(2n)$ 
%(\resp\ $I_1\oplus I_1$ if $n=1$; see \cref{fig:2d})
(the $4n$-gons), $B_d$ and $F_4$ are \emph{not} uniquely determined up to scale and orientation, but each case forms a continuous 1-dimensional family of combinatorially equivalent zonotopes (see \eg\ \cref{fig:B3_mod} for the case $B_3$).
%
%\begin{myenumerate}
%\item
%the $4n$-gons form a continuous 1-dimensional family of combinatorially~equi\-valent vertex-transitive zonotopes,
%\item
%the zonotopes of the $B_d$-family form a continuous 1-dimensional family of combinatorially equivalent vertex-tran\-sitive zonotopes, and
%\item
%the $F_4$-zonotopes form a continuous 1-dimensional family of combinatorially equivalent vertex-tran\-sitive zonotopes.
%\end{myenumerate}
%
%The mentioned cases are special for the following reason: among the root systems, only the vectors in $I_2(2n)$ (or $I_1\oplus I_1$ if $n=1$; see \cref{fig:2d}), $B_d$ and $F_4$ do not form a single, but \emph{two} orbits under their symmetry group.
%The reason for the difference lies in the action of the Weyl group on the corresponding root system.
%The reason for these exceptions is as follows: 
Typically, the vectors of a root system form a single orbit under the action of the associated Weyl group, except in the cases $I_2(2n)$, $B_d$ and $F_4$ in which they form two orbits \cite[Section 2.11]{humphreys1990reflection}.
%
%The mentioned cases are special, as these are these correspond to the only root systems for which the con
%The reason for this difference is that the 
The length of the vectors in each orbit can be chosen independently from the other orbit, giving each such zonotope one degree of freedom that manifests in two (possibly) different edge lengths.

Such degrees of freedom are also present in all reducible vertex-transitive zonotopes.
For example, the $d$-cubes (the cartesian product of line segments) belong to the continuous family of $d$-orthotopes with $d$ degrees of freedom.
%Equivalently, hyperprisms are reducible, as they are the cartesian product of polygons and line segments.

%All vertex-transitive zonotopes of the same combinatorial type have the same normalization (\cf\ \cref{def:homogeneous}).
%These normalizations are homogeneous, and are exactly the polytopes known as \emph{omnitruncated uniform polytopes} \cite{...}.
%The enumeration of those amounts to listing all Coxeter-Dynkin diagrams in which \emph{all} nodes are ringed.

%\subsubsection{Two dimensions}
%
%There exist infinitely many 1-parameter families of vertex-transitive zonotopes in two dimensions.
%
%For each $n\ge 2$, the regular $2n$-gon is a vertex-transitive zonotope.
%This belongs to the reflection groups $I_2(n)$ if $n\ge 3$, or $I_1\oplus I_1$ if $n=2$.
%
%If $n$ is odd, this is the only vertex-transitive zonotope with $2n$ vertices.
%But if $n$ is even, we actually have an infinite 1-dimensional family of vertex-transitive $2n$-gonal zonotopes.
%The reason is, that then $\Gen(Z)$ has two orbits of vectors \wrt\ $\Aut(\Gen(Z))$, and their lengths can be chosen independently.
%This results in two possibly different edge lengths in the zonotope (see \cref{fig:2d}).

We give further information on some of the families:

\subsubsection{The family $A_d$}
%\subsubsection{Families that exist in all dimensions}
%In sufficiently high dimensions there are exactly three types of (irreducible) vertex-transitive zonotopes resulting from the reflection groups $A_d, B_d$ and $D_d$.
%In small dimensions, some of these might coincide, or further exceptional ones might exist.
%Among these three classes, only the $B_d$-permutahedron provides an infinite 1-dimensional family.

The generated zonotope is the standard permutahedron.
%This vertex-transitive zonotope is uniquely determined up to scale and orientation.
For $d=3$, this zonotope is called \emph{truncated octahedron} (\cref{fig:permutahedron}), which coincides with the zonotope obtained from $D_3$ (this is the only pair of coinciding zonotopes).
The classical name for this family is \emph{omnitruncated $d$-simplices}.

\subsubsection{The family $B_d$}
%The root system of $B_d$ has two orbits under its symmetry group, which corresponds to two (possibly) different edge lengths in the corresponding zonotope.
The vertices of a general $B_d$-permutahedron are formed by the coordinate permutations and sign selections of some vector
\begin{equation}\label{eq:Bd_coordinates}
(\pm x_1,\..., \pm x_d)\in\RR^d,\qquad \text{with $x_1,\...,x_d>0$}.
\end{equation}
If the $x_i$ form a linear sequence $x_i=x_0 + \eps (i-1)$ for some fixed $x_0,\eps>0$, then the corresponding polytope is a zonotope.
The quotient $\eps/x_0$ parametrizes the 1-dimensional family (see \cref{fig:B3_mod}), and $\eps/x_0=\sqrt 2$ corresponds to the homogeneous representative.
Clasically, the homogeneous zonotopes of this family are called \emph{omnitruncated $d$-cubes}.

\subsubsection{The family $D_d$}
%The $D_d$-root system is formed by one of the orbits of the $B_d$-root system ($B_d$ contains $2d^2$ vectors and decomposes into two orbits of size $2d^2-2d$ and $2d$).
%The zonotope generated by the remaining vectors is a vertex-transitive zonotope and is uniquely determined up to scale and orientation.
The vertex-coordinates of the $D_d$-zonotope are as in \eqref{eq:Bd_coordinates} with setting $x_0=0$ in the linear sequence (see \cref{fig:B3_mod}).
For $d=3$, this zonotope coincides with the $A_3$-permutahedron (the truncated octahedron).
For $d=4$, this zonotope is known as the \emph{truncated 24-cell}.
Classically, the zonotopes of this family are called \emph{omnitruncated $d$-demicubes}.

\subsubsection{Exceptional zonotopes}

In dimensions $d\in\{3,4,6,7,8\}$ the following exceptional vertex-transitive zonotopes exist:

\begin{myenumerate}
\item 
For $d=3$, there exists the \emph{omnitruncated icosahedron/dodecahedron} to the reflection group $H_3$ (see \cref{fig:permutahedron}).
\item 
For $d=4$, there exists the continuous 1-dimensional family of $F_4$-zonotopes (the homogeneous member is called \emph{omnitruncated 24-cell}), as well as the \emph{omnitruncated 120-cell/600-cell} to the reflection group $H_4$.
\item
For $d\in\{6,7,8\}$, there exist the omnitrucation of the uniform $E_d$-polytope.
\end{myenumerate}

\subsubsection{Summary}

All in all, we obtain the following numbers of irreducible homogeneous zonotopes (aka.\ combinatorial types of irreducible vertex-transitive zonotopes) per dimension.

\vspace{0.3em}
\begin{center}
\begin{tabular}{r||c|c|c|c|c|c|c|c|c|l}
$d$ & 0 & 1 & 2        & 3 & 4 & 5 & 6 & 7 & 8 & $\ge 9$  \\
\hline
\# &  1 & 1 & \clap{$\infty$} & 3 & 5 & 3 & 4 & 4 & 4 & 3
\end{tabular}
\end{center}
\vspace{0.3em}

\section{Characterizing root systems}
\label{sec:root_systems}

Our results on vertex-transitive and homogeneous zonotopes enable us to give interesting alternative characterizations of root systems:

%\Cref{res:zonotope_equivalents} provides us with an alternative characterization of root systems:
%Let $R\subset\RR^d$ be a finite set of vectors. 
%A \emph{half} of $R$ is the subset $R(c):=\{r\in R\mid \<c,r\> > 0\}\subseteq R$, for which $c$ is non-orthogonal to all $r\in R$.

\begin{corollary}\label{res:root_halves}
If $R\subset\RR^d\setminus\{0\}$ is finite and reduced, then the following are equivalent:
\begin{myenumerate}
	\item $R$ is a root system, and
	\item all semi-stars of $R$ are congruent.
\end{myenumerate}
\begin{proof}
%We assume that $R$ is full-dimensional, otherwise we restrict to $\Span(R)$.
$R$ is reduced, hence the set of generator of $\Zon(R)$.
We can apply $(ii)\Leftrightarrow (iii)$ from \cref{res:zonotope_equivalents}.
%
%$(i)$ is equivalent to $\Zon(i)$ being a $\Gamma$-permutahedron, and $(ii)$ is equivalent to $\Zon(R)$ being vertex-transitive.
%It is easy to see that 
%If $R$ is reduced, the proof follows immediately from \cref{res:zonotope_equivalents} (from the equivalence $(ii)\Leftrightarrow (iii)$).
%The non-reduced case follows. \TODO
%\msays{What if $R$ is not centrally symmetric?}
%
%If $(i)$ or $(ii)$, then $R$ is centrally symmetric.
%This is clear for $(i)$.
%To see this for $(ii)$ choose a generator $r\in R$ with $-r\not\in R$, and a vector $c\in\RR^d$ orthogonal to $r$, but to no other vectors in $R\setminus\Span\{r\}$. 
%For a sufficiently small $\eps>0$ the halves $R_{c\pm \eps r}$ are identical with the difference that one contains $r$, and the other one does not.
%In particular, they have different numbers of elements and cannot be congruent.
%
\end{proof}
\end{corollary}

The statement of \cref{res:root_halves} still holds true, even if $R$ is not reduced, but we leave that to the reader.

The \emph{norm} of a semi-star shall be the norm of the sum of its elements.
It follows from \cref{res:root_halves}, that a root system has all semi-stars of the same norm.
We prove a weaker form of a converse:

\begin{theorem}\label{res:root_norms}
Let $R\subset \mathbb S^d$ be a finite centrally symmetric set of unit vectors.
If all semi-stars of $R$ have the same norm, then $R$ is a root system.
%\msays{Do we need central symmetry?}
%
\begin{proof}
Recall that each vertex of $\Zon(R)$ can be written as the sum of the vectors in some semi-star of $R$ (\cref{res:semi_star_vertices}). 
The norm of that semi-star therefore equals the distance of that vertex from the origin.

Since all semi-star have the same norm, $\Zon(R)$ has all vertices on a common sphere around the origin.
Furthermore, $R$ is centrally symmetric, and all vectors in $R$ are of the same length, \ie\ all edges of $\Zon(R)$ are of the same length.
Thus, $\Zon(R)$ is homogeneous, and by \cref{res:homogeneous_root_system} it is a $\Gamma$-permutahedron.

Since $R$ is centrally symmetric and consists of unit vectors, it is reduced.
Conclusively, $R=\Gen(\Zon(R))$ is a set of generators of a $\Gamma$-permutahedron, and therefore a root system.
\end{proof}
\end{theorem}

It is unclear to the author whether dropping the assumption of central symmetry in \cref{res:root_norms} allows other, essentially different vector configurations (see \cref{q:unit_vectors_root_system}). %, or whether it follows from all semi-stars having the same norm (see \cref{q:unit_vectors_root_system}).
On the other hand,~condition $R\subset\Sph^{d}$ is certainly necessary, as there are inscribed zonotopes with edges of distinct lengths and that are not permutahedra (see \cref{fig:inscribed_zono} and \cref{sec:open}).

%% file: sec-homo/conclusion.tex
\section{Related problems and open questions}
\label{sec:open}

%\begin{figure}
%\centering
%\includegraphics[width=0.45\textwidth]{img/comb_permutahedron}
%\caption{A zonotope that has all vertices on a common sphere, but is not a $\Gamma$-permutahedron. The corresponding generators have all semi-stars of the same norm, but form no root system.}
%\label{fig:comb_permutahedron}
%\end{figure}

\begin{figure}
\centering
\includegraphics[width=0.8\textwidth]{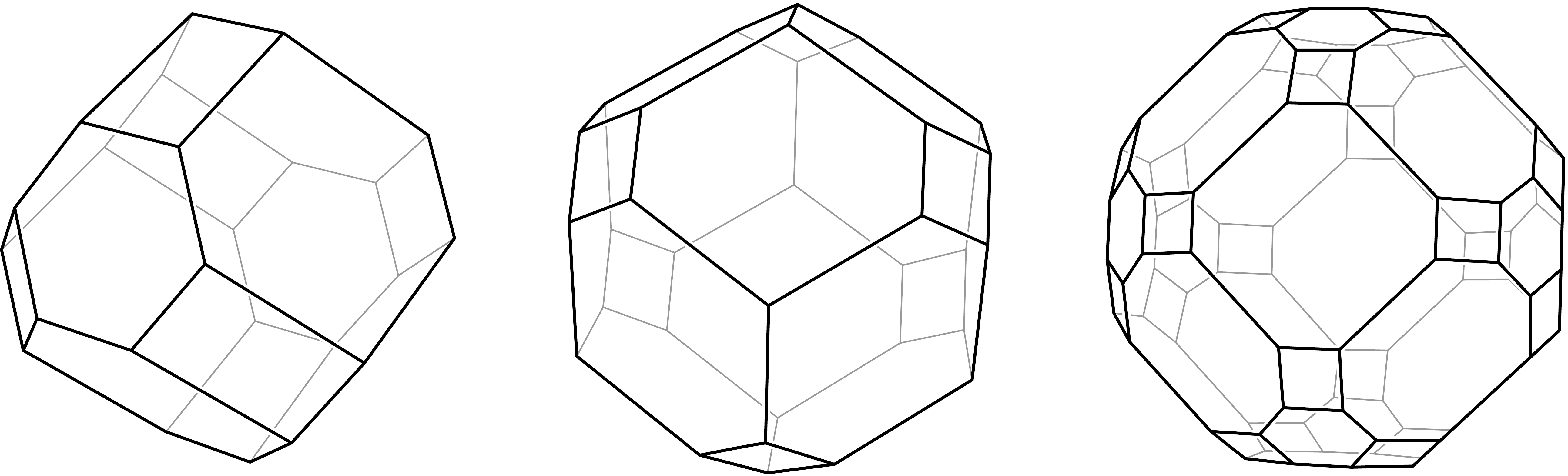}
\caption{Three inscribed zonotopes obtained as projections of higher-dimensional $\Gamma$-permutahedra. Left: projection of the $A_4$-zonotope, combinatorially equivalent to the $A_3$-zonotope, but \emph{not} vertex-transitive. Middle: projection of the $D_4$-zonotope. Right: projection of the homogeneous $F_4$-zonotope.
The latter two are not combinatorially equivalent to a $\Gamma$-permutahedron.}
\label{fig:inscribed_zono}
\end{figure}

The vertex-transitive zonotopes belong to the larger class of \emph{inscribed zonotopes} (that is, all vertices are on a common sphere).
If we put no restrictions on the edges (as in the homogeneous case), we obtain a class with much unclearer properties.

It is clear that not all inscribed zonotopes can be $\Gamma$-permutahedra: it is easy enough find an inscribed centrally symmetric $2n$-gon that is not vertex transitive.
It is harder to construct examples in three or more dimensions.
One might be temped to conjecture that an inscribed zonotope is at least \emph{combinatorially equivalent} to a $\Gamma$-permutahedron.
However, counterexamples were provided by Raman Sanyal and Sebastian Manecke (personal communication): the orthogonal projection of a $\Gamma$-permutahedron along one of its edge directions is again inscribed, but not necessarily combinatorially equivalent to a $\Gamma$-permutahedron (see \cref{fig:inscribed_zono}).
In general, the projection of an inscribed zonotope along an edge direction is again inscribed.
There are further known examples which cannot be obtained as such repeated projections of $\Gamma$-permutahedra, but all those are still combinatorially equivalent to one which was obtained as a projection.

\begin{question}[by R. Sanyal and S. Manecke]
Are there inscribed zonotopes which are not combinatorially equivalent to a $\Gamma$-permutahedron or the repeated projection of such one along edge directions?
\end{question}

Zonotopes have a known relation to real hyperplane arrangement.
The results of this paper immediately yields the following: a hyperplane arrangement whose symmetry group acts transitively on its chambers must be a reflection arrangement.
A more general question was asked by Caroline J. Klivans and Ed Swartz \cite[Problem 13]{klivans2011projection}, and we shall repeat it here:

\begin{question}[by C. Klivans and E. Swartz] \label{q:hyperplane_arrangement}
If all chambers of a real hyperplane arrangement are congruent, is it a reflection arrangment?
\end{question}

One finds that such an arrangement must be central and simplicial.
The answer  to \cref{q:hyperplane_arrangement} is known to be affirmative in dimensions $d\in\{2,3\}$ \cite{ehrenborg2016coxeter}, but is open in $d\ge 4$.
Dualizing again, the analogous question for zonotopes is the following:

\begin{question}
If all vertex-figures of a zonotope are identical, is it combinatorially equivalent to a $\Gamma$-permutahedron, or more precisely, is its normalization (see \cref{def:normalization}) a $\Gamma$-permutahedron?
\end{question}

%One reasonable question might be the following:
%
%\begin{question}
%If a zonotope has all vertices on a common sphere, is it combinatorially equivalent to a $\Gamma$-permutahedron?
%\end{question}
%
%The answer is trivially positive in two dimensions. 
%See \cref{fig:comb_permutahedron} for a zonotope in which all vertices are on a common sphere, but which is \emph{not} a $\Gamma$-permutahedron.
%The answer is not clear in $\ge 3$ dimensions.

Another question, for which no immediate answer was found, is whether it is necessary to assume central symmetry in \cref{res:root_norms}, or whether the following stronger version of the theorem holds:

\begin{question}\label{q:unit_vectors_root_system}
Let $S\subset\mathbb S^d$ be a finite set of unit vectors, in which all semi-stars have the same norm. Is $S\cup -S$ necessarily a root system?
\end{question}

Note that $S$ alone is not necessarily a root system: let $S$ be the set of vertices of a regular triangle centered at the origin, then all semi-stars have the same norm, but $S$ is not a root system.
However,  $S\cup -S$ is the root system of $I_2(3)$.

%% file: sec-homo/appendix.tex
\appendix

\section{}
\label{sec:appendix}

%\stepcounter{section}

We provide proofs of \cref{res:semi_star_vertices} and \cref{res:linearly_closed_faces}.
These can be found in the standard literature on zonotopes (e.g. \cite{mcmullen1971zonotopes}), but are included for completeness.

We need the following preliminary result:

\begin{proposition}\label{res:gen_faces}
Given a face $\sigma\in\F(Z)$, there is a unique partition of $\Gen(Z)=R_-\cupdot R_0\cupdot R_+$ into disjoint sets, so that
\begin{equation}\label{eq:face_R}
\sigma=\sum_{\mathclap{r\in R_+}} r + \Zon(R_0).
\end{equation}
If $\sigma$ is the set of maximizers of $\<\free,c\>$ in $Z$ for $c\in\RR^d$, then
$$R_0:=\{r\in \Gen(Z)\mid \<r,c\>=0\},\quad R_{\pm}:=\{r\in\Gen(Z)\mid \pm \<r,c\>>0\},$$
and these sets are independent of the exact choice of $c$.
\begin{proof}
By \eqref{eq:sum_of_generators}, each point $v\in Z$ can be written in the form
$$v=\sum_{\mathclap{r\in \Gen(Z)}} a_r r,\quad \text{with coefficients $a_r\in[0,1]$}.$$
Fix some arbitrary partition $\Gen(Z)=R_-\cupdot R_0\cupdot R_+$. Consider the following two separate sets of constraints for the coefficients $a_r$ of $v$:
$$
(*)\; a_r=\begin{cases}
\phantom+0 & \text{for $\<r,c\> < 0$}\\
\phantom+1 & \text{for $\<r,c\> > 0$}\\
\text{arbitrary} & \text{for $\<r,c\> = 0$}
\end{cases},
\quad
(**)\;a_r=\begin{cases}
\phantom+0 & \text{for $r\in R_-$}\\
\phantom+1 & \text{for $r\in R_+$}\\
\text{arbitrary} & \text{for $r\in R_0$}
\end{cases}.
$$
%
%Let $\sigma^{(*)}$ \resp\ $\sigma^{(**)}$ be the set of all points in $Z$ which satisfy $(*)$ \resp\ $(**)$.
The constraints in $(*)$ are equivalent to $v$ maximizing the functional $\<\free, c\>$ in $Z$.
On the other hand, the constraints in $(**)$ are equivalent to
$$
v \in 
\Big\{\sum_{\mathclap{r\in R_+}} r + \sum_{\mathclap{r\in R_0}} a_r r\,\Big\vert\, a\in[0,1]^{R_0}\Big\}
=\sum_{\mathclap{r\in R_+}}r + \Zon(R_0).$$
By definition, the face $\sigma$ is the set of all points in $Z$ that satisfy $(*)$.
If $\sigma$ is of the form \eqref{eq:face_R}, then $\sigma$ is additionally the set of all points in $Z$ that satisfy $(**)$.
This happens exactly when $(*)$ and $(**)$ are equivalent conditions.
Comparing the definitions of $(*)$ and $(**)$, wee see that we necessarily need to choose $R_-,R_0$ and $R_+$ as claimed.
\end{proof}
\end{proposition}

\renewcommand{\thesection}{\arabic{section}}
\setcounter{section}{2}
\setcounter{theorem}{2}

\begin{proposition}%\label{res:semi_star_vertices}
The vertices $\F_0(Z)$ are in one-to-one correspondence with the semi-stars $S\subset\Gen(Z)$ via the bijection
$S\; \mapsto\;  v_S:=\sum_{r\in S}r.$
\begin{proof}
By \cref{res:gen_faces}, for each vertex $v\in \F_0(Z)$ there are \emph{uniquely} determined disjoint subsets $R_0,R_+\subseteq\Gen(Z)$ with
$$v=\sum_{\mathclap{r\in R_+}} r + \Zon(R_0)$$
Clearly $R_0=\eset$, and $R_+$ is then a semi-star.
\end{proof}
\end{proposition}

\begin{proposition}%\label{res:linearly_closed_faces}
The following holds:
\begin{myenumerate}
\item 
For each face $\sigma\in\F(Z)$, $\Gen(\sigma)\subseteq\Gen(Z)$ is a flat.
\item
For each flat $F\subseteq\Gen(Z)$ there is a face $\sigma\in\F(Z)$ with $\Gen(\sigma)=F$.
\end{myenumerate}
\begin{proof}
Let $\sigma\in \F(Z)$ be some face. 
By \cref{res:gen_faces}, $\sigma$ is a translate~of $\Zon(R_0)$ for some subset $R_0\subseteq\Gen(Z)$. In particular, $R_0$ is a flat.
Then $(i)$ follows via
$$
(*)\quad \Gen(\sigma)=\Gen(\Zon(R_0))=R_0.
$$
%
%The first equivalence uses that translation does not change the generators.
%The second equivalence uses that $R_0$ is reduced. 

For $(ii)$, choose a vector $c\in\RR^d\nozero$, so that $F=\Gen(Z)\cap c^\bot$, which is possible since $F\subseteq\Gen(Z)$ is a flat.
The set $\sigma$ of maximizers of $\<\free,c\>$ in $Z$ is a face of $Z$.
We apply \cref{res:gen_faces} to $\sigma$ to obtain the set $R_0\subseteq\Gen(Z)$, which turns out to be exactly $F$.
Then $\Gen(\sigma)=R_0=F$ follows via $(*)$.
\end{proof}
\end{proposition}